\documentclass{article}
\usepackage{graphicx}
\usepackage{amsmath,amsthm,amssymb}
\usepackage{xcolor}
\newif\ifFigs \Figstrue
\Figsfalse

\newcommand{\Includegraphics}[2]%
{\centering
    \ifFigs\includegraphics[width=#1]{#2}%
    \else  \includegraphics[draft,width=#1]{#2} \fi}

\newcommand{\includesubgraphics}[3]%
{\begin{minipage}{#1}
        \includegraphics[width=\textwidth]{#2}\\
        \centerline{{\footnotesize (#3)}}
    \end{minipage}
}

\usepackage[normalem]{ulem}

\definecolor{blue-violet}{rgb}{0.54,0.17,0.89}
\definecolor{amethyst}{rgb}{0.6,0.4,0.8}
\definecolor{darkviolet}{rgb}{0.58, 0.0, 0.83}
\definecolor{darkgreen}{rgb}{0,.4,0}
\definecolor{mixedgreen}{rgb}{0.3,0.6,00}
\definecolor{bananayellow}{rgb}{1.0, 0.88, 0.21}
\definecolor{arylideyellow}{rgb}{0.91, 0.84, 0.42}
\definecolor{bananamania}{rgb}{0.98, 0.91, 0.71}

\newcommand{\REMOVEva}[1]%
           {{\color{magenta}\sout{#1}}}

\begin{document}

\begin{center}

{\LARGE\bf Bistability in a one-dimensional model of a
  two-predators-one-prey population dynamics system}

\bigskip

{\bf Sergey Kryzhevich$^{1,2,3*}$, Viktor Avrutin $^4$, Gunnar S\"oderbacka$^{5}$}   

\bigskip

{\small

$^{1}$ \quad Institute of Applied Mathematics, Faculty of Applied Physics and Mathematics, Gdańsk University of Technology, ul. Narutowicza 11/12, 80-233 Gda\'{n}sk, Poland\\

$^{2}$ BioTechMed Center, Gda\'{n}sk University of Technology, ul. Narutowicza 11/12, 80-233 Gda\'{n}sk, Poland\\

$^{3}$ \quad Saint Petersburg State University,
14 line of the VO, house 29B; Saint Petersburg, 199178 Russia\\

$^{4}$ \quad Institute for Systems Theory and Automatic Control, University of Stuttgart, Pfaffenwaldring 9, 70569 Stuttgart, Germany; viktor.avrutin@ist.uni-stuttgart.de\\

$^{5}$ \quad \r{A}bo Akademi,
Turku FI-20500, Finland; gsoderba@abo.fi\\

$^{*}$ Corresponding author: kryzhevicz@gmail.com}
\end{center}	

{\small 
\noindent\textbf{Abstract.} 
In this paper, we study the classical two-predators-one-prey model. The classical model described by a system of 3 ordinary differential equations can be reduced to a one-dimensional bimodal map. We prove that this map has at most two stable periodic orbits. Besides, we describe the structure of bifurcations of the map. Finally, we describe a mechanism that yields bistable regimes. We find several areas of bistability numerically.

\bigskip

}

\section{Introduction}

Modeling ecosystems is one of the well-established parts of the theory of dynamical systems. One of the standard models in this area is the Lotka-Volterra system, also known as the prey-predator model. The behavior of this two-dimensional autonomous system of differential equations is well known: the system may exhibit stationary or periodic solutions only.

Conversely, ecological models involving competition between two or more predators (and, possibly, more species of prey) may demonstrate a sophisticated behavior, including chaotic dynamics. Besides, the number of attracting sets (periodic or chaotic) may vary.

The general case of a two-predators-one-prey system is considered by Hsu, Hubell, and Waltman in \cite{Hsu1} -- \cite{Hsu3}.

The models considered in the cited works are based on ecological laws described initially by Holling in \cite{holling}. A good review of recent results in population dynamics can be found in \cite{diekmann}.

A study of periodic solutions of the above-mentioned models together with a detailed bifurcation analysis is presented in \cite{bw}, \cite{keener} and \cite{smith}.  Two-dimensional and three-dimensional models in discrete time and the bifurcations of their equilibria are considered in \cite{wk}.

A model with additional terms responsible for competition between predators is studied in \cite{sswa}. Extinction conditions are discussed for predator species, stability of fixed points and other invariant sets is analyzed.

A model with ratio-depending predator growth rates is introduced in \cite{du}. Equilibrium analysis is performed, and extinction conditions for predators are presented.

Various generalizations of the Lotka-Volterra model with polynomials at the right-hand side have also been studied, see, for example, \cite{hmb} and \cite{sg}.

More sophisticated models including stochastic terms \cite{fdbh}, \cite{lm} or taking into account diffusion phenomena \cite{ml} are also developed.

In the present work, we consider a classical two-predators-one-prey model. Following \cite{sp}, this model can be approximated by a one-dimensional bimodal map. Previously, it has been observed that typically this map has a unique attractor, although for some parameter values it also may exhibit bistability (coexistence of two attractors). The goal of the present paper is to describe the mechanism leading to this effect.

The paper is organized as follows.

In Sec.~\ref{sec:model} and Sec.~\ref{sec:map}, the considered 3D model in continuous time and its one-dimensional approximation in discrete time are introduced, respectively. For the one-dimensional map, we apply the classical techniques developed by Devaney and demonstrate that this map cannot have more than two attractors. Thereafter, in Sec.~\ref{sec:bistability} we identify the regions in the parameter space where bistability occurs. We demonstrate that this kind of dynamics occurs in neighborhoods of the intersection points of bifurcation curves forming so-called shrimp-structures.  
\section{\label{sec:model}Two-predators-one-prey model} 

Let us consider the two-predators-one-prey model given by the following three differential equations
\begin{equation}\label{eq1}
\begin{array}{l}
  \dot X_i=p_i \varphi_i(S) X_i-d_iX_i, \qquad i=1,2;\\
  \dot S=H(S)-q_1 \varphi_1(S) X_1-q_2 \varphi_2(S) X_2
\end{array}
\end{equation}
investigated previously in \cite{sp}. Here, the non-negative values $S$ and $X_i$, $i=1,2$, represent quantities of the prey and predators respectively; $H$ and $\varphi_i$ are smooth functions with $H(0)=0$, $H$ being of the logistic type and both $\varphi_i$ being non-decreasing. Parameters $p_i$, $d_i$ and $q_i$ are positive.

Using the normalized variables
\begin{equation*}
    s=\frac{S}{K},\quad y_i=\frac{q_iX_i}{rK},\qquad i=1,2;
  \end{equation*}
Eqs.~\eqref{eq1} can be rewritten as
\begin{equation}\label{eq2}
\begin{array}{l}
\dot y_i=\phi_i(s)y_i, \qquad i=1,2;\\
\dot s= h(s) - \psi_1(s)y_1-\psi_2(s)y_2 
\end{array}
\end{equation}
where 
$$h(s)=\dfrac{1}{K} H(sK), \qquad \psi_i(s)=\varphi_i(sK), \qquad \phi_i(s)=p_i \psi_i(s)-d_i.$$ 

In this paper, we choose the functions
$\phi_i(s)$, $h(s)$, and 
$\psi_i(s)$ to be defined as follows:
  \begin{equation*}
    \begin{split}
      \phi_i(s)&=m_i\frac{s-\lambda_i}{s+a_i},\qquad i=1,2;\\
      \psi_i(s)&=\frac{s}{s+a_i},\qquad i=1,2;\\
      h(s) &=(1-s)s
    \end{split}
  \end{equation*}
where $a_1$, $a_2$, $\lambda_1$ and $\lambda_2$ are positive. Then, system \eqref{eq2} takes the following form: 
\begin{equation}\label{eq3}
\begin{array}{l}
\dot y_i=m_i\dfrac{s-\lambda_i}{s+a_i}y_i, \qquad i=1,2\\
\dot s= \left(1-s - \dfrac{y_1}{s+a_1}-\dfrac{y_2}{s+a_2}\right)s. 
\end{array}
\end{equation}

Dissipativity of this system and the extinction conditions for one of the predators are considered in \cite{sp} (see also \cite{malaysians}). There, possible equilibria of this system are studied as well as some periodic solutions. In particular, it is shown in the cited work that there are no predator coexistence stationary solutions for
$\lambda_1\neq \lambda_2$ or $a_1\neq a_2$ and that the coexistence of periodic or chaotic solutions is possible.

The structure of Poincar\'{e} maps corresponding to the condition $s = \mathrm{const}$, $\dot s < 0$ near fixed points is studied in
\cite{osipovIJBC}.
It is shown in \cite{eirola1} and \cite{eirola2} that for a broad range of parameter values the considered system exhibits a strong contraction in the $(y_1 + y_2)$-direction in which case its dynamics can be approximated by the one-dimensional map given by
\begin{equation}
x_{n+1}=f(x_n) = \beta + x_n -\dfrac{k_1+k_2 e^{x_n}}{1+e^{x_n}}u.
\label{initialone-dimensional}
\end{equation}
where $\beta$,$u$ and $k_i$ are parameters. Here $x_j = \log(y_{2j}/y_{1j})$ where $y_{1j}$ and $y_{2j}$ are values of $y_1$ and $y_2$ respectively calculated at the $j$-th step.

\bigskip

\noindent\textbf{Remark 1.} Note that for parameter values which do not lead to a strong contraction in the $(y_1 +y_2)$-direction, system \eqref{eq3} cannot be described by a one-dimensional map and presumably may exhibit such phenomena as Shil'nikov's spiral chaos and coexistence of at least three attractors (see \cite{bkk,osipovIJBC,sp}). As discussed below, this cannot occur in the one-dimensional
model \eqref{initialone-dimensional}.

\section{\label{sec:map}One-dimensional model and its properties}

\begin{figure}[t!]
\begin{center}
\includegraphics[height=2in]{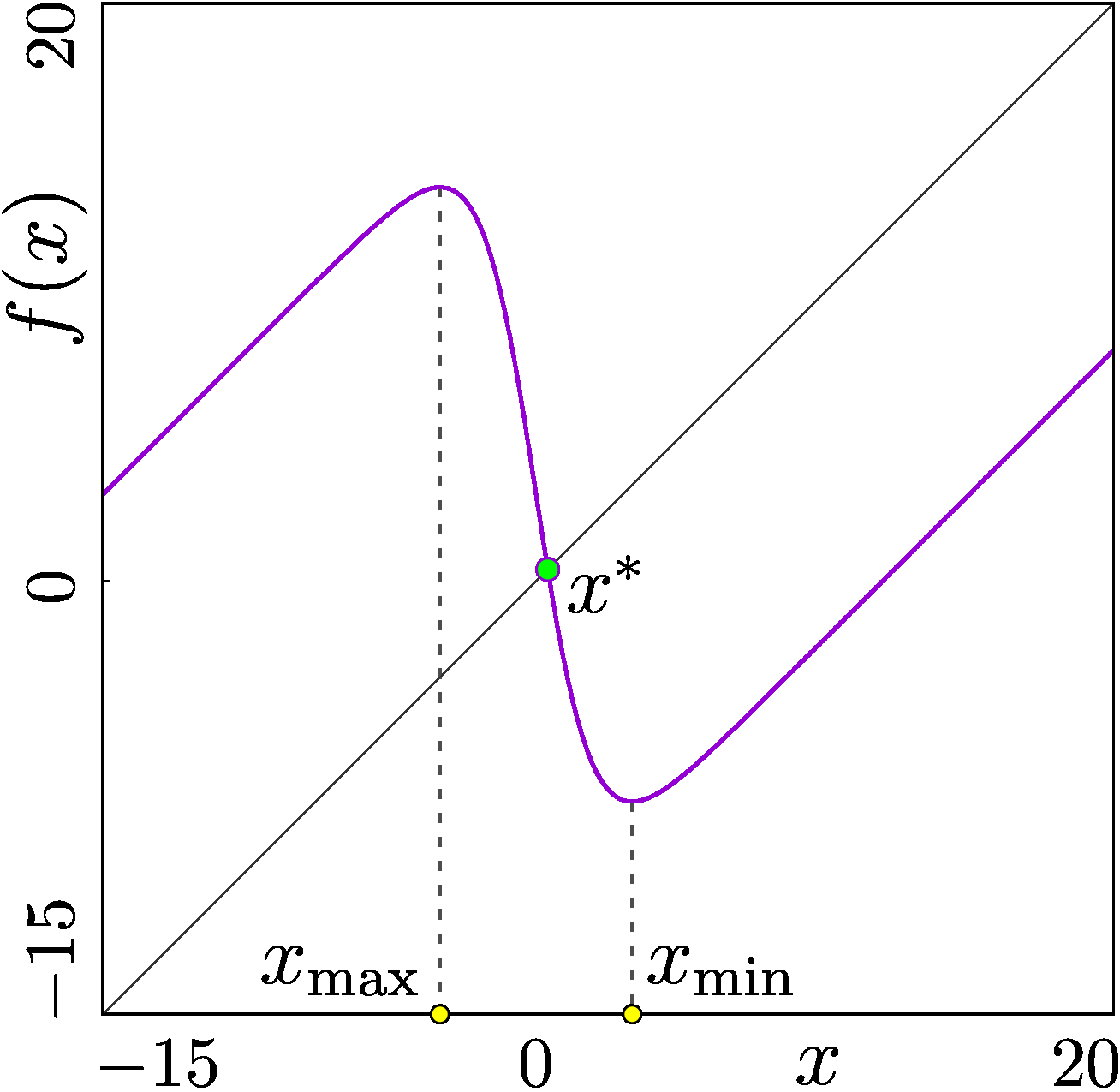}
\caption{\footnotesize Graph of the one-dimensional map~\eqref{glavUrav}, $b=-12$, $k=-30$.}
\end{center}
\label{fig:graph}
\end{figure}

It can easily be shown that map~\eqref{initialone-dimensional} can be written in the following form
\begin{equation}
x_{n+1}=f(x_n) = b + x_n -\frac{k}{1+e^{x_n}}
\label{glavUrav}
\end{equation}
where $b=\beta+k_2 u$, $k=k_1u$. An example of the graph of the function $f$ in Eq.~\eqref{glavUrav} is shown in Fig.~\ref{fig:graph}. In the following we consider the dynamics of map~\eqref{glavUrav} in the parameter domain
\begin{equation}
  \label{eq:domain:P}
    P=\{ (b,k) \mid k<b<0\}
\end{equation}
as outside this domain all trajectories either diverge or converge to a stable fixed point.

\begin{figure}[t!]
\begin{center}
  \includegraphics[width=.84\textwidth]{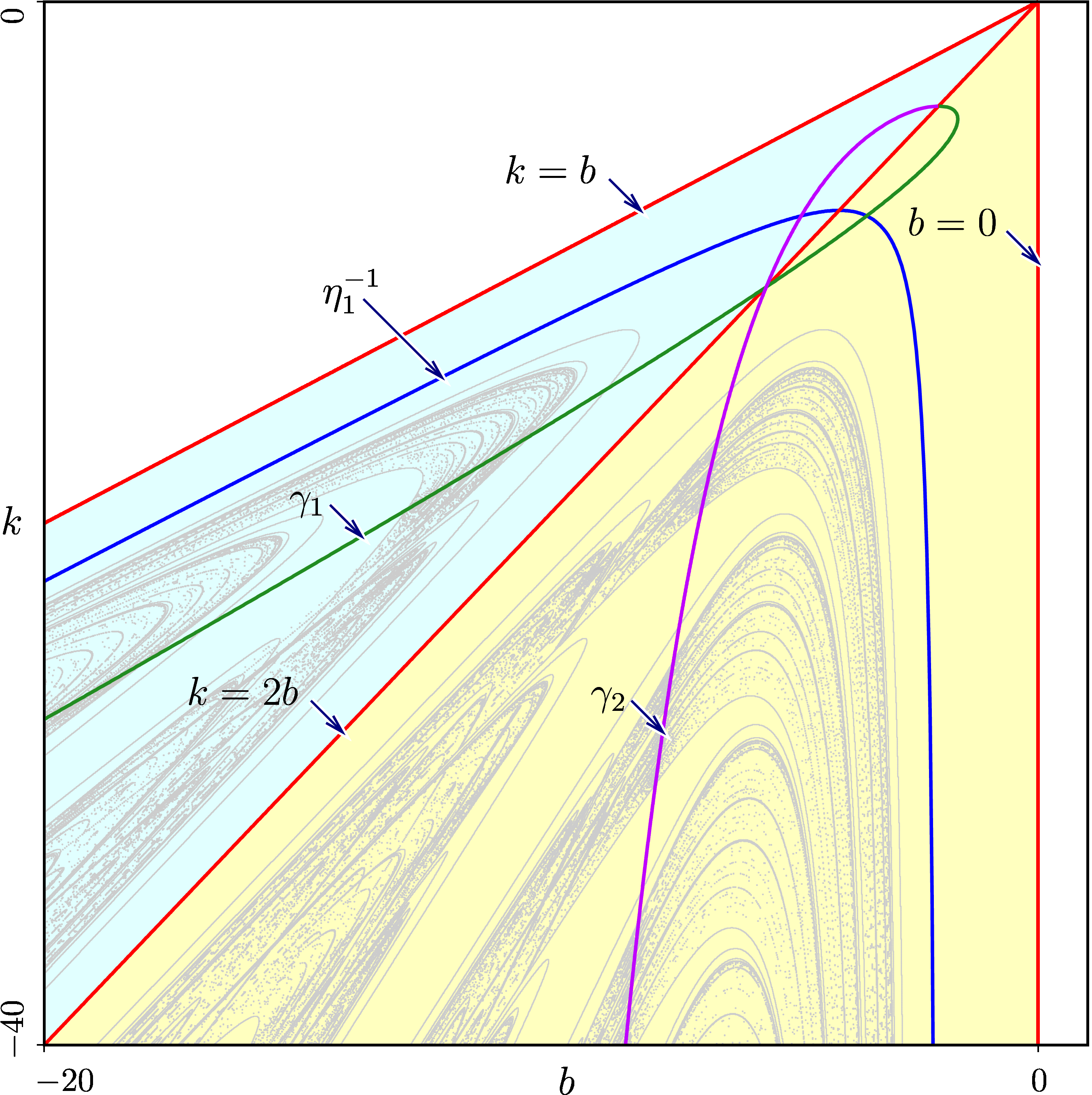}
\end{center}
\caption{\footnotesize Boundaries of the regions $P$, $P_1$, $P_2$, $P^s_{x°*}$ in the parameter space and the boundaries of the regions corresponding to different definitions of the globally attracting absorbing intervals. Additionally, the bifurcation structures calculated numerically is shown (see also Fig.~\ref{bkregpp})}
\label{fig:curves}
\end{figure}

\noindent\textbf{Remark 2.} The map has the following symmetry property:
$$\left.f(x)\right|_{b=b^*, k=k^*} =
  \left.-f(-x)\right|_{b=k^*-b^*, k=k^*}$$
Therefore, if the map at the parameter values $b=b^*$, $k=k^*$ has a period-$n$ orbit $\{x_1,\ldots,x_n\}$, then the map at the parameter values $b=k^*-b^*$, $k=k^*$ has the period-$n$ orbit $\{-x_1,\ldots,-x_n\}$. As a consequence, the bifurcation structures of the map in the parameter regions
$$P_1 = \{(b,k) \mid 2\, b < k < b < 0\} \quad
\text{and} \quad P_2 = \{(b,k) \mid k\leq 2\, b < 0 \}$$
with $P = P_1 \cup P_2$ (see Fig.~\ref{fig:curves}) are topologically equivalent.

\noindent\textbf{1. The fixed point and the period-2 orbit.} 

Notice that for $k\ge -4$ the map is 
For $k<-4$, the function $f$ has a local minimum and a local maximum at the points
  \begin{equation*}
    \begin{split}
      x_{\max}&=\ln\left(-1-\tfrac{1}{2}k - \tfrac{1}{2}\sqrt{4k+k^2}\right),\\
      x_{\min}&=\ln\left(-1-\tfrac{1}{2}k + \tfrac{1}{2}\sqrt{4k+k^2}\right)
    \end{split}
  \end{equation*}
respectively. Evidently, $x_{\min}=-x_{\max}>0$. The map is increasing for $x<x_{\max}$ and for $x>x_{\min}$ and decreasing in the interval $\lbrack x_{\max},\, x_{\min}\rbrack$ containing the unique fixed point
\begin{equation*}
  x^*=\ln \left( \frac{k}{b} -1 \right)
\end{equation*}
existing for $(b,k)\in P$.  As the parameter values approach the boundaries of $P$, the fixed point tends to $\pm \infty$ (i.e., for a fixed value of $k<0$, we have $x^*\to -\infty$ if $b \to k-0$ and $x^* \to \infty$ if $b\to 0$.
  
From the condition $f'(x^*)>-1$, we obtain that the fixed point is attracting (moreover, it is globally attracting) in the parameter region
\begin{equation*}
  P^{s}_{x^*} = \left\{(b,k)\in P \mid k>\frac{b^2}{b+2}, \quad b>-2 \right\}.
\end{equation*}
At the boundary of the region $P^{s}_{x^*}$, i.e., at the curve
\begin{equation*}
  \eta^{-1}_1 = \left\{ (b,k)\in P \mid k=\frac{b^2}{b+2}, \quad b>-2\right\}
\end{equation*}
the fixed point $x^*$ undergoes a supercritical period-doubling bifurcation leading to the appearance of a period-2 orbit. This can be checked by direct calculations. As shown by the following lemma, this period-2 orbit exists in the complete parameter region $P \setminus P^{s}_{x^*}$:

\bigskip

{\noindent\textbf{Lemma 1.}} \emph{For any $(b,k)$ such that $k<\frac{b^2}{b+2}$, $b<-2$, the map \eqref{glavUrav} has a period-2 orbit}.
    
\bigskip

\noindent\textbf{Proof.} Evidently, 
$$\lim_{x\to-\infty} f(x)-x=b-k, \qquad 
\lim_{x\to\infty} f(x)-x=b.$$ 
Consequently,
$$\lim_{x\to-\infty} f^2(x)-x=2(b-k)>0, \qquad 
\lim_{x\to\infty} f^2(x)-x=2b<0.$$
On the other hand, $(f^2)'(x^*)=(f'(x^*))^2>1$. Therefore, the function $f^2(x)-x$ has at least one zero in the interval $(-\infty,x^*)$ and another one in the interval $(x^*,\infty)$. The fixed point $x^*$ is unique, so these zeros correspond to a period-2 orbit.
$\square$

Moreover, one can actually prove a more general result:

\noindent\textbf{Lemma 2.} \emph{For any $(b,k)$ such that $k<\frac{b^2}{b+2}$, $b<-2$, at least one of the following statements applies:
    \begin{enumerate}
    \item There exists an $n\in {\mathbb N}$ such that map \eqref{glavUrav} has a stable period-$2^n$ orbit.
    \item The map has period-$2^n$ orbits for all $n\in {\mathbb N}$.
    \end{enumerate}}

The proof of this lemma is similar to that of Lemma 1.

Lemma 1 proves the existence of at least one period-2 orbit. However, for each $b \in (P \setminus P^{s}_{x^*})$ there exists only one period-2 orbit. This is shown in Lemma 3.

\bigskip

\noindent\textbf{Lemma 3.} \emph{For any $(b,k)$ such that $k<\frac{b^2}{b+2}$, $b<-2$, map \eqref{glavUrav} can have at most one period-2 orbit.}

\bigskip

\noindent\textbf{Proof.}  The idea of the proof is as follows. Given a period-2 orbit 
$$\lbrace x_1,\, x_2\rbrace , \, x_2 >x_1$$ of map~\eqref{glavUrav}, we express $k$ and $b$ as functions of $x_1$ and $x_2$. Then we consider values $k$ and $b$ such that 
$$B:= \frac{2\, b}{k}=\mathrm{const},$$ 
and represent $k$ as a function of $x_1$. We show that the latter function is strictly monotone which implies that two distinct orbits $\{x_1,x_2\}$ cannot correspond to a pair $(k,b)$.

Now let us provide detailed proof. Solving equations $f(x_1)=x_2$ and $f(x_2)=x_1$ for $b$ and $k$ we get
\begin{equation}
  \begin{split}
B =& u_1 + u_2, \,\, k=\frac{2\, (x_2 -x_1)}{u_2-u_1},\\
&\text{where}\quad u_i=\frac{1}{1+e^{x_i}},\quad i=1,\, 2. 
  \end{split}
\label{bk}
\end{equation}
Note that $0<u_i<1$ and $u_2<u_1$ implying $u_1>B/2$. From $k<b<0$ follows $0<B<2$. Thus $B/2 < u_1 < B$ if $B<1$ and $B/2 < u_1 < 1$ if $B>1$. Solving the last equality of \eqref{bk} for $x_i$ we get $x_i = \ln \left( \frac{1}{u_i} -1\right) = h(u_i)$. Henceforth, we drop the subscript and write $u$ instead of $u_1$. So, we can write down $k$ as a function of $u$:
$$k(u)= \frac{2(h(u)-h(B-u))}{2u-B}.$$
Since 
$$h'(u)=\frac{1}{u\, (u-1)}$$ 
we obtain 
\begin{equation}\label{incond}
\lim_{u\to B/2} k(u) = 2\, h'(B/2) =\frac{8}{B\, (B-2)},
\end{equation}
where the limit is taken from the right side. Moreover, it follows by the Lagrange theorem that
$k< 2 \max h'(u) =  h'(B/2) =\frac{8}{B\, (B-2)}.$
The function $h'(u)$ is convex for $u\in (0,1)$ and, consequently 
$$h'(B-u)+h'(u)<2h'(B/2)=\dfrac{8}{B\, (B-2)}$$
for any $u\in (B/2,\min(B,1))$.

Now we calculate the derivative of $k$:
$$k'(u)= \frac{2}{2u-B}  \left(h'(B-u)+h'(u) - k \right)< \frac{2}{B-2u}  \left(\frac{8}{B\, (B-2)} - k \right)$$
for any $u>B/2$. So, taking into account initial conditions \eqref{incond}, we obtain that
$$k(u)<\frac{8}{B\, (B-2)}$$ 
and, moreover, $k'_u(u)<0$ for any $u>B/2$. This implies the uniqueness of the period-2 orbit. $\square$

\noindent\textbf{Remark 3.} Orbits (even stable ones) of periods higher than 2 may be non-unique. One can find numerically that two period-4 attractors can coexist.

\begin{figure}[t!]
\begin{center}
  \includegraphics[width=.84\textwidth]{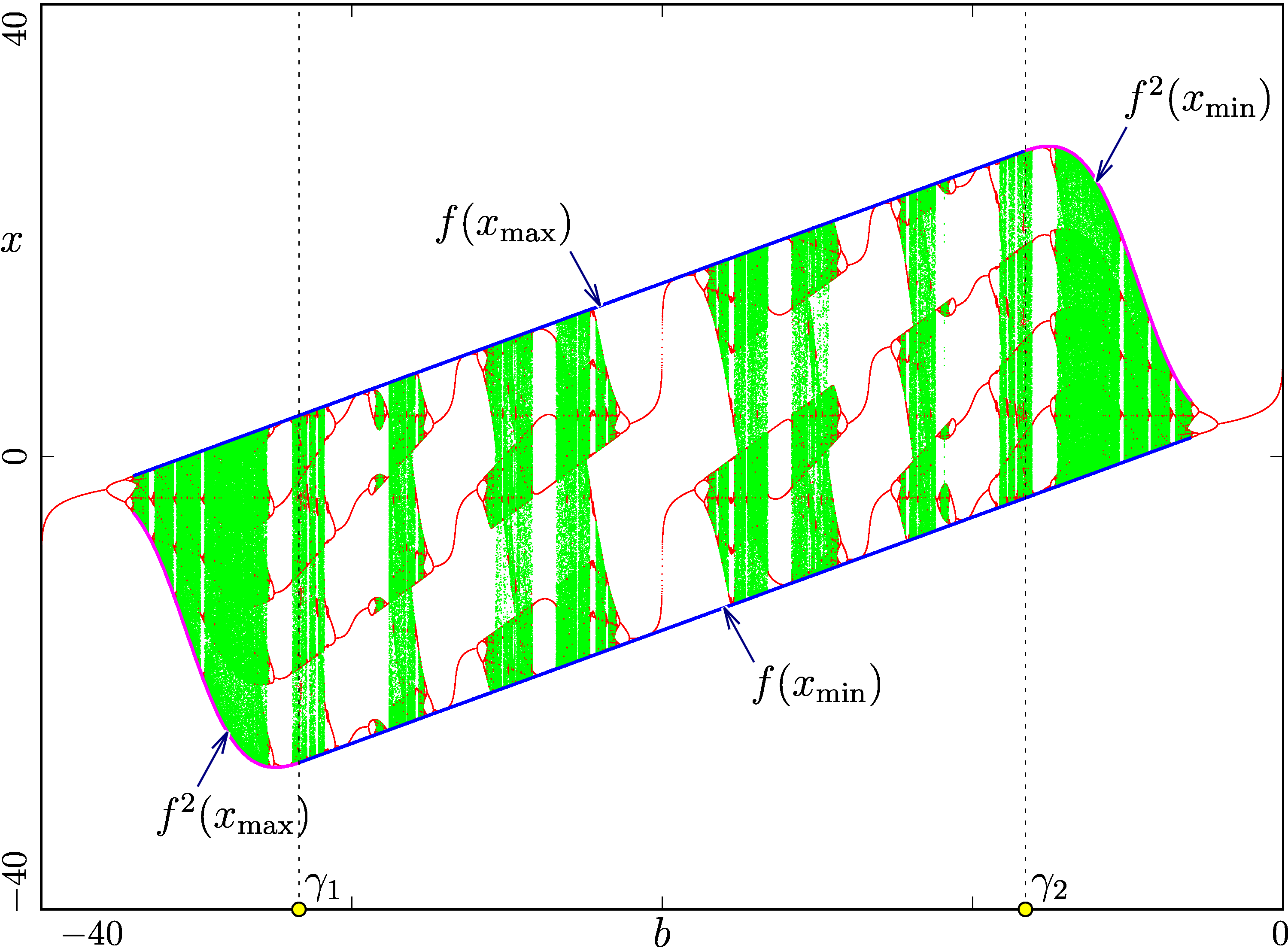}
\end{center}
\caption{\footnotesize Bifurcation sequence at $k=-40$. Boundaries of the
  globally attracting absorbing intervals are indicated.}
\label{fig:bif:one-dimensional}
\end{figure}

\noindent\textbf{2. A globally attracting interval.}  It is easy to see that for $(b,k)\in P$ the function $f$ satisfies $f(x)>x$ for $x<x^*$ and $f(x)<x$ for $x>x^*$. Therefore, either the fixed point is globally attracting or there map has a globally attracting absorbing interval around $x^*$. The boundaries of this interval are given by the images of the points $x_{\max}$ and $x_{\min}$, as illustrated in Fig.~\ref{fig:bif:one-dimensional}.  As one can see, in the left part of this figure, the absorbing interval is given by $[f^2(x_{\max}), f(x_{\max})]$, in the middle part by $[f(x_{\min}), f(x_{\max})]$, and in the right part by $[f(x_{\min}), f^2(x_{\min})]$.  As shown below in Lemma 4, the regions in the parameter space corresponding to these configurations are separated from each other by the curves
    \begin{align*}
    \begin{split}
      \gamma_1 =& \left\{ (b,k)\in P \mid f(x_{{\max}}) = x_{\min} \right\}
      = \left\{ (b,k)\in P \mid b=b_1\right\}\\
      &\text{where} \quad
     b_1= 2\, x_{\min} -1 - e^{x_{\min}}\\
      &\phantom{\text{where} \quad b_2}
       =2\ln\left(-1-\tfrac{1}{2}k + \tfrac{1}{2}\sqrt{4k+k^2}\right)+\tfrac{1}{2}k - \tfrac{1}{2}\sqrt{4k+k^2};  
    \end{split}\\
    \begin{split}
      \gamma_2 =& \left\{ (b,k)\in P \mid f(x_{\min})>x_{{\max}} \right\}
      =\left\{ (b,k)\in P \mid b=b_2\right\}\\
      &\text{where} \quad
    b_2= -2\, x_{\min} -1 - e^{-x_{\min}}\\
      &\phantom{\text{where} \quad b_1}
       =2\ln\left(-1-\tfrac{1}{2}k - \tfrac{1}{2}\sqrt{4k+k^2}\right)+\tfrac{1}{2}k + \tfrac{1}{2}\sqrt{4k+k^2}.
    \end{split}
    \end{align*}
{(see Fig.~\ref{fig:curves} for graphs of $\gamma_1$ and $\gamma_2$)}. Using this notation, we can state the following:

\bigskip

\noindent\textbf{Lemma 4.} \emph{
\begin{enumerate}
\item If $b<\min (b_1,b_2) $ then the interval $J_-=\lbrack (f^2(x_{\max}),\, f(x_{\max})\rbrack$ does not contain $x_{\min}$ and is globally attracting.
\item If $b>\max (b_1,b_2)$ then the interval $J_+=\lbrack (f(x_{\min}),\, f^2(x_{\min})\rbrack$ does not contain $x_{\max}$ and is globally attracting.
\item If $b_2<b < b_1$ then the interval $\lbrack (f(x_{\min}),\, f(x_{\max})\rbrack$ does not contain any of $x_{\min}$ and $x_{\max}$. In this case, the function $f$ is monotonous on this interval. Therefore, the function can have a fixed point and a period-2 orbit only. 
\item If $b_1<b < b_2$ then the interval $J_0=\lbrack (f(x_{\min}),\, f(x_{\max})\rbrack$ contains both $x_{\min}$ and $x_{\max}$ and is globally attracting.
\end{enumerate}} 

\bigskip

The proof of this lemma follows from the definitions of $b_1$ and $b_2$.

\bigskip

Note that the absorbing intervals $J_-$, $J_+$ and $J_0$ in cases 1, 2 and 4 respectively are positively invariant.

\begin{figure}[t!]
\begin{center}
  \includegraphics[width=.84\textwidth]{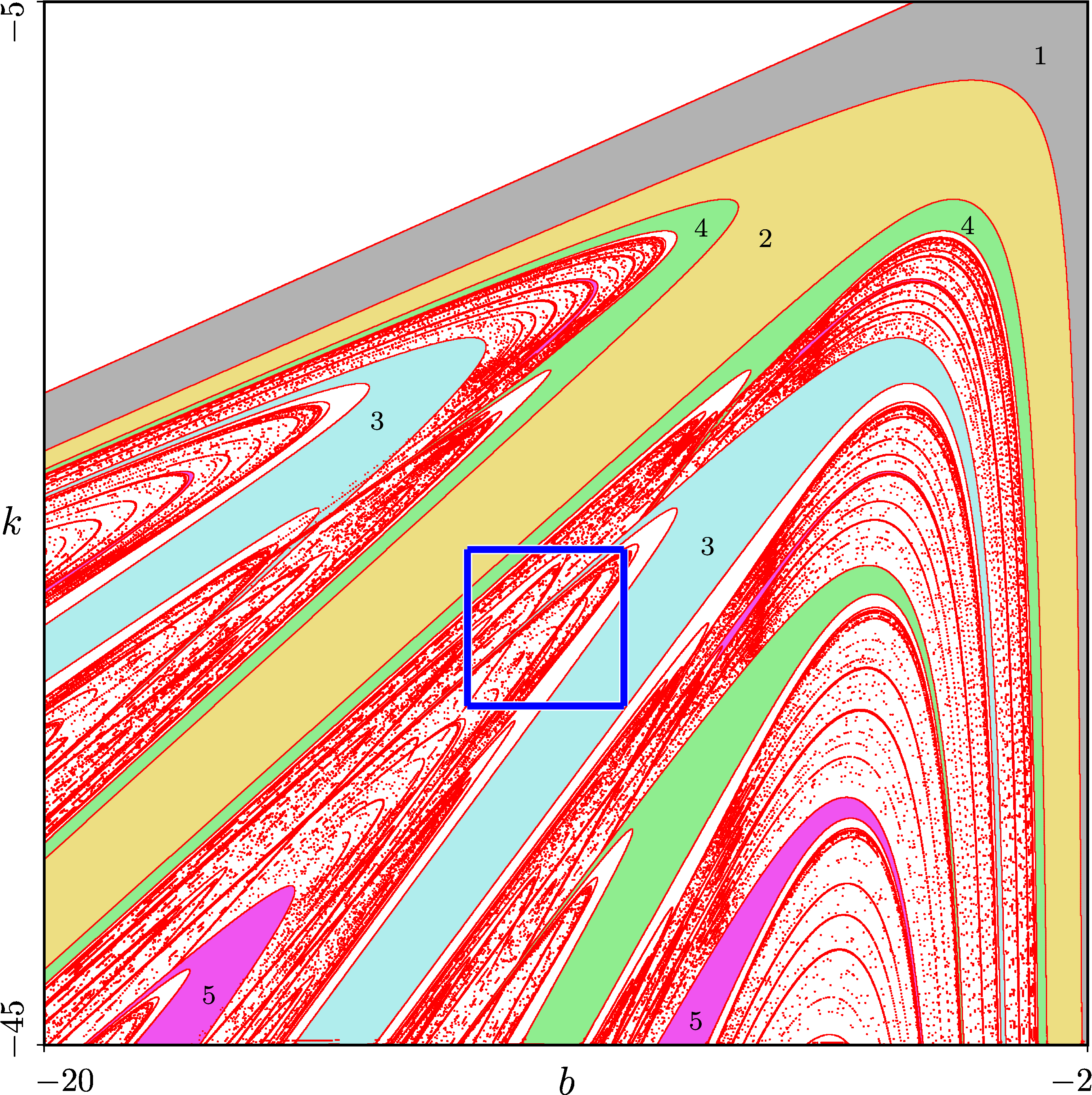}
\end{center}
\caption{\footnotesize Regions of the parameter plane corresponding to different periods of stable orbits.}
\label{bkregpp}
\end{figure}

However, if $b_2 <b <b_1$, the fixed point and the period-2 orbit cannot be stable at the same parameter value. Hence, in this case, there is only one stable periodic orbit. 

Calculating the Schwarzian derivative for the map \eqref{glavUrav}, we get
$$
S \, f(x)= k \, e^x \, \frac{2(e^x-1)^2-(k+4)\, e^x}{2((e^x+1)^2+ k\, e^x)^2}
$$
\noindent which is negative for $k<-4$. 

Let $b<\min (b_1,b_2)$. Then the positively invariant and globally attracting absorbing interval 
$J_-=[f^2(x_{\max}),f(x_{\max})]$ contains the point $x_{\max}$ and does not contain $x_{\min}$. In other words, the map is unimodal on $J_-$. It follows from the classical result of Devaney \cite[Theorem 11.4]{devaney} that in this case any stable periodic orbit (and, in fact, any attractor) attracts $x_{\min}$. So, such an orbit is necessarily unique. The case $b\ge\max(b_1,b_2)$ is similar. Therefore, the dynamics of the map for parameter values inside
the region $P$ but outside the region between the curves $\gamma_1$, $\gamma_2$ below their intersection point (see Fig. \ref{fig:curves}) cannot be affected by bistability.

Finally, if $b_1 <b <b_2$, the globally attracting positively invariant absorbing interval $J_0$ contains both $x_{\min}$ and $x_{\max}$. Therefore, as follows from the Devaney's result mentioned above, it can contain two attractors, with $x_{\min}$ belonging to the basin of one of them, and $x_{\min}$ to the basin of other one. Accordingly, the map in the parameter region shown in Fig. \ref{fig:curves} between the curves $\gamma_1$, $\gamma_2$ below their intersection point can exhibit bistability.

\section{\label{sec:bistability}Period doubling cascades and coexisting \\ periodic solutions of distinct periods}

\begin{figure}[t!]
\begin{center}
  \includesubgraphics{.4\textwidth}{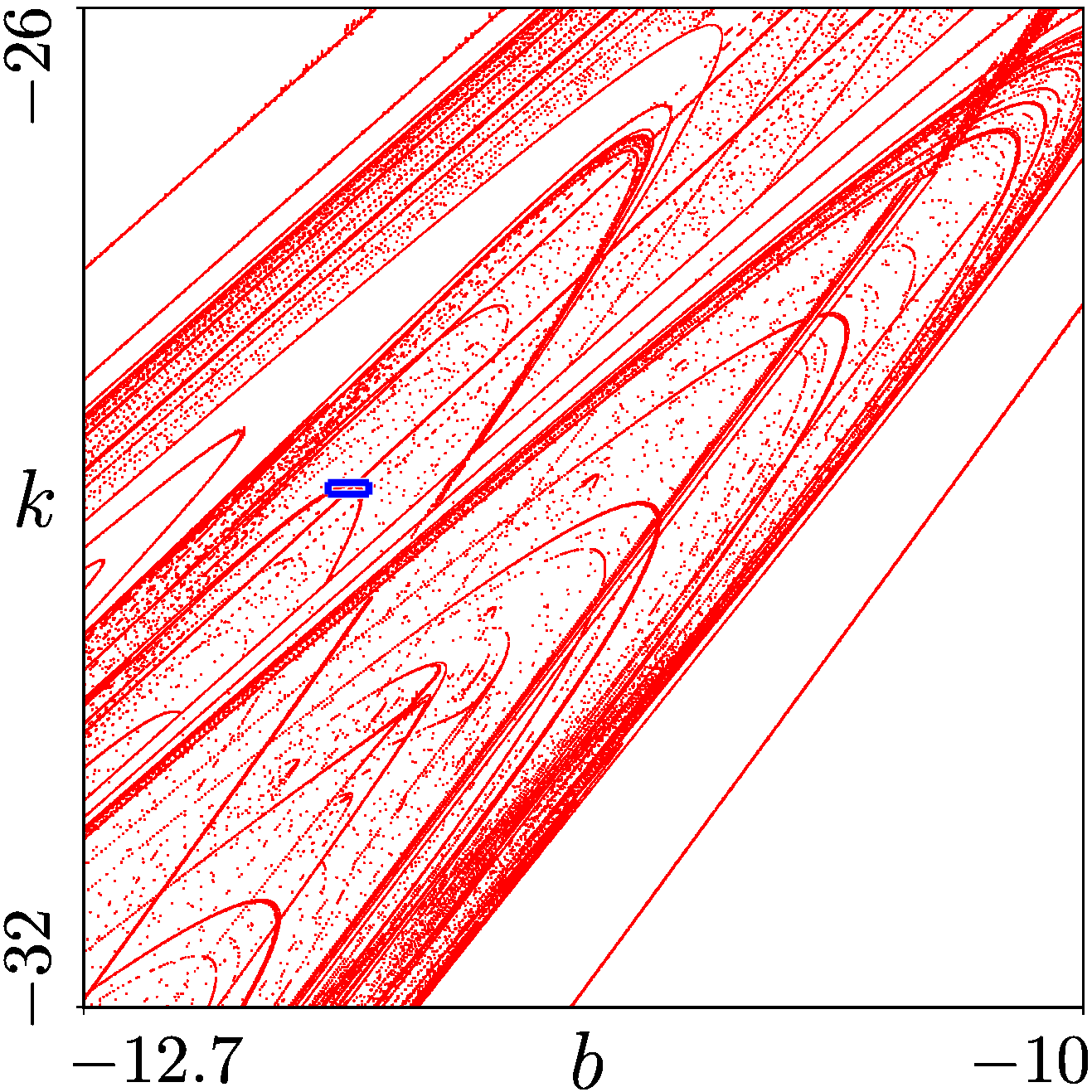}{a}\quad
  \includesubgraphics{.4\textwidth}{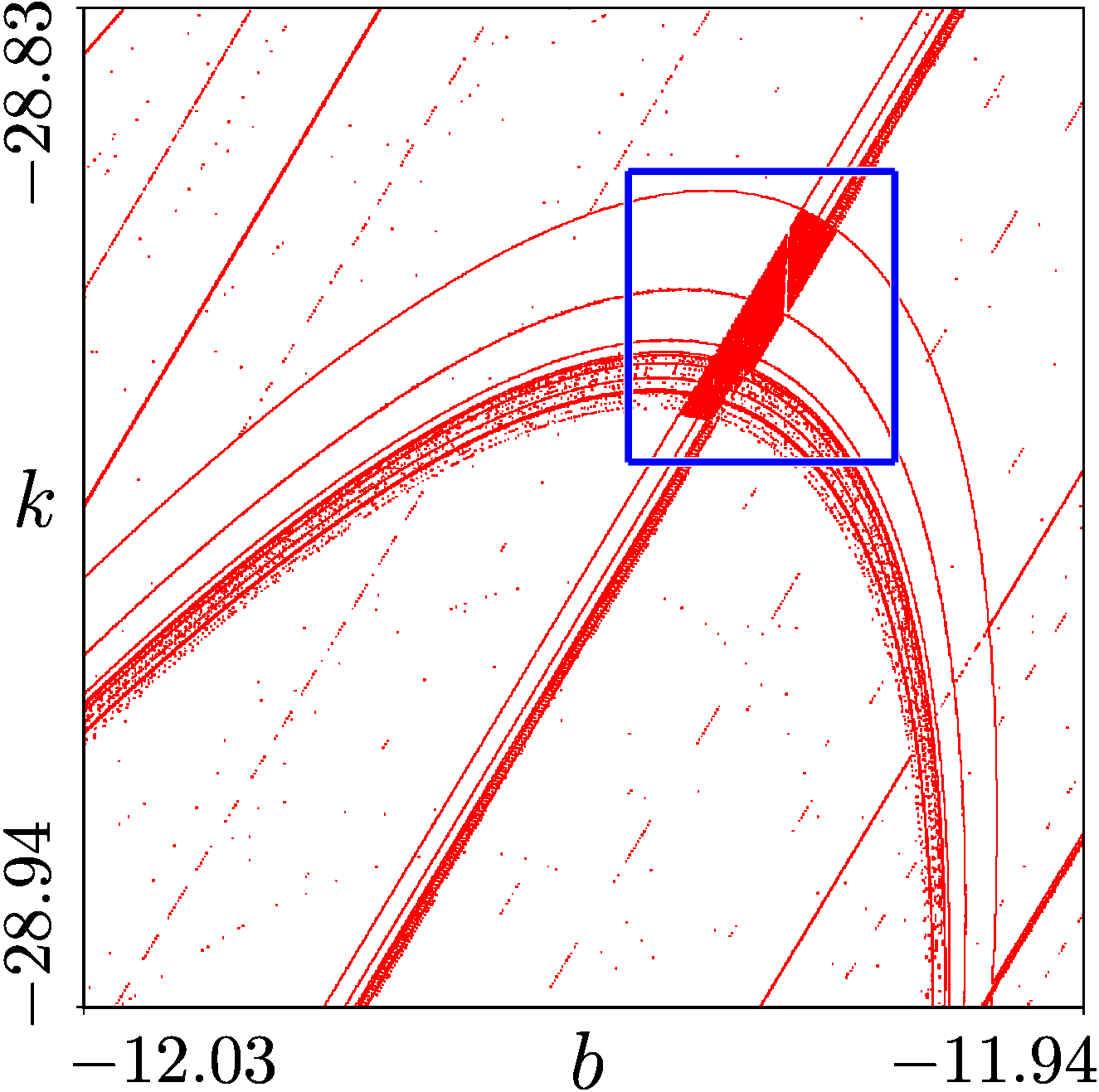}{b}\\[1ex]
  \includesubgraphics{.84\textwidth}{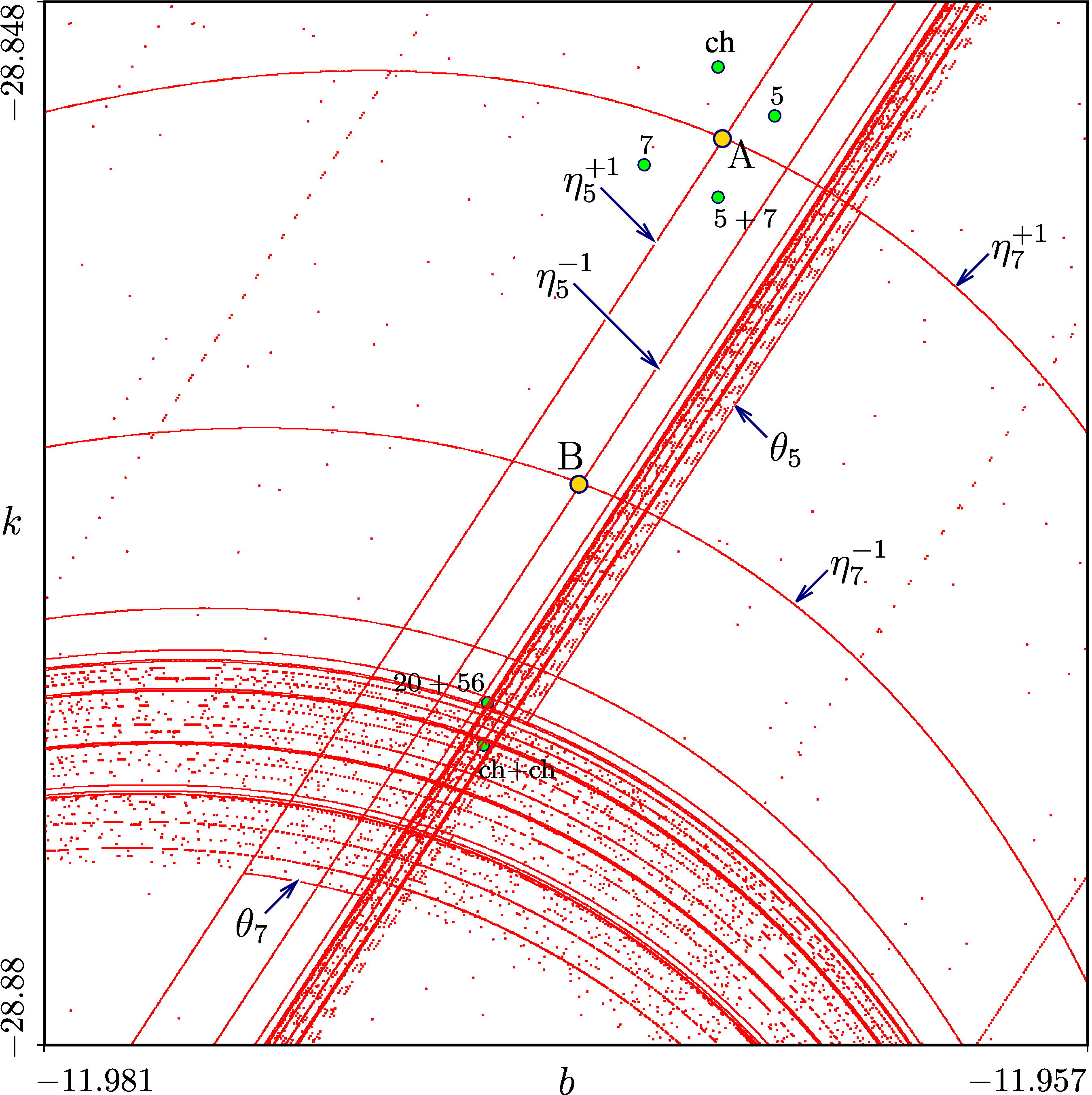}{c}
\end{center}
\vspace*{-1.5em}
\caption{\footnotesize Bifurcation structure in the $(b,k)$ parameter plane.  (a)
  Magnification of the rectangle marked in Fig.~\ref{bkregpp}.  (b)
  Magnification of the rectangle marked in (a).  (c) Magnification of the rectangle marked in (b).  Fold bifurcation curves $\eta^{+1}_5$, $\eta^{+1}_7$, flip bifurcation curves $\eta^{-1}_5$, $\eta^{-1}_7$, and final bifurcation curves $\theta_5$, $\theta_7$ are indicated. Attractors at the marked parameter points are shown in Fig.~\ref{fig:attractors:A} and Fig.~\ref{fig:attractors:C}.  }
\label{fig:BxK:blowup}
\vspace*{-1em}
\end{figure}

In order to explain the occurrence of bistability in
map~\eqref{glavUrav} in the parameter region between the curves $\gamma_1$, $\gamma_2$, let us consider the complete bifurcation structure in the 2D parameter space $(b,k)$. Fig.~\ref{bkregpp} shows this structure including regions corresponding to stable cycles of higher periods calculated numerically. As illustrated in the magnification of this structure shown in Fig.~\ref{fig:BxK:blowup}, in the parameter domain $P$ given by Eq.~\eqref{eq:domain:P} the bifurcation curves form so-called \emph{shrimp}-structures~\cite{Gallas1994}, previously
observed in several one-dimensional and two-dimensional
maps~\cite{Oliveira2011,Stoop2012}, including the well-known H\'enon map. For a detailed description of the ``anatomy'' of a shrimp-structure, we refer to~\cite{Glass2013}.

A distinguishing feature of the shrimp-structures is the ``tails'', i.e., long and narrow parameter regions confined from one side by a fold bifurcation curve issuing from the central part of the structure and approaching infinity. Inside these narrow regions, one observes a complete period-doubling cascade and thereafter by the complete logistic map scenario, including a countable set of curves associated with quasi-periodic dynamics (Feigenbaum-attractors) as well as by an uncountable set of curves related to non-robust chaotic dynamics. From the other side, the ``tails'', are confined by expansion or final bifurcations (see, e.g., \cite{VA-BOOK2019}) of narrow-band chaotic attractors (interior and boundary crises, respectively) associated with homoclinic bifurcations of the unstable cycles appearing at the fold bifurcations.

It is well-known that such ``tails'' may overlap pairwise.  In this case, two transversely intersecting fold bifurcation curves subdivide the parameter plane into four quadrants. In one of these quadrants, the attractors belonging to both overlapping ``tails'' coexist pairwise, which explains the occurrence of bistability in map~\eqref{glavUrav}. Accordingly, the region of bistability related to an intersection of two ``tails'' is confined by four bifurcation curves (two-fold bifurcations and two homoclinic bifurcations).

\begin{figure}[t!]
\begin{center}
  \includesubgraphics{.4\textwidth}{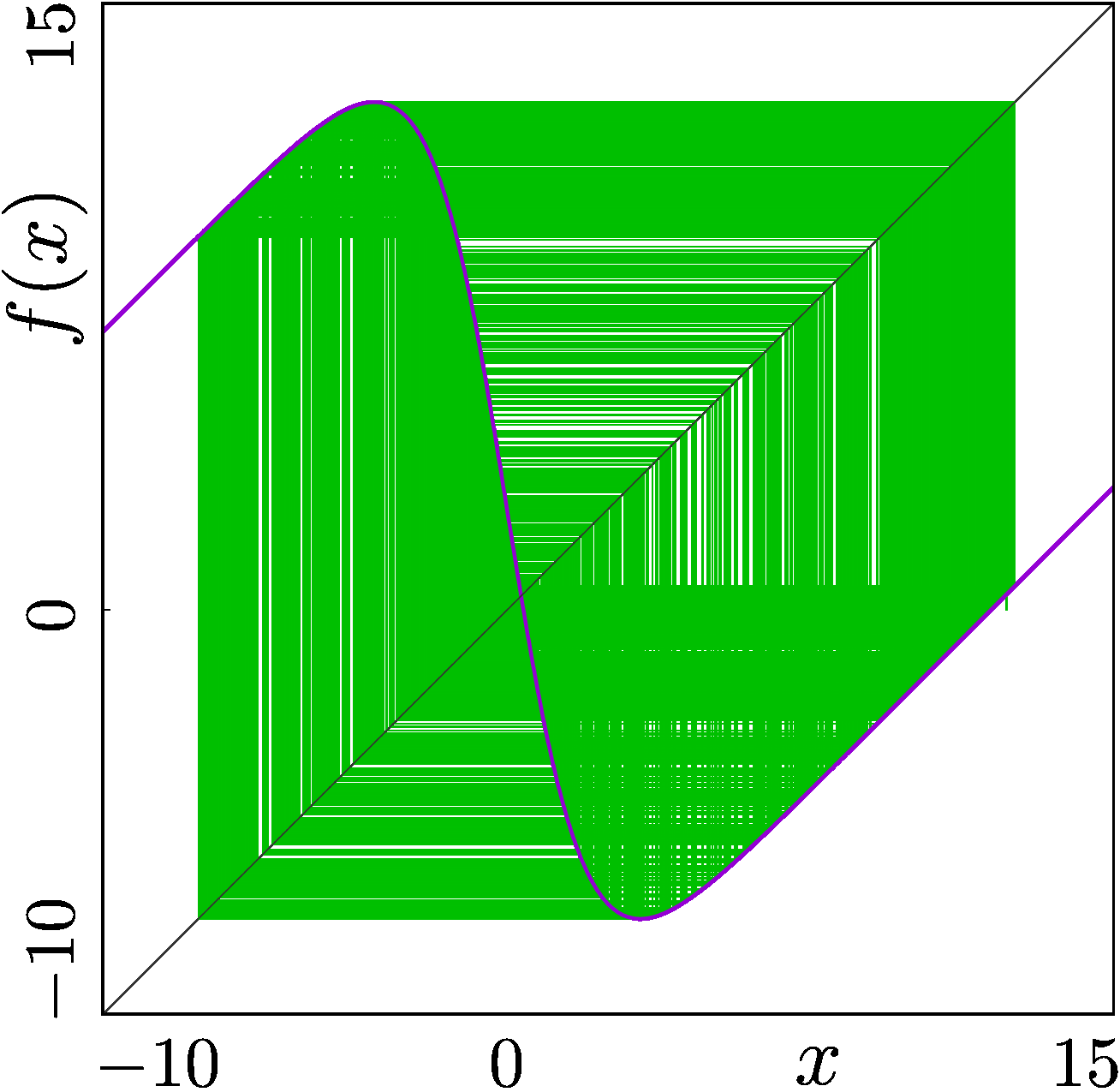}{a}\quad
  \includesubgraphics{.4\textwidth}{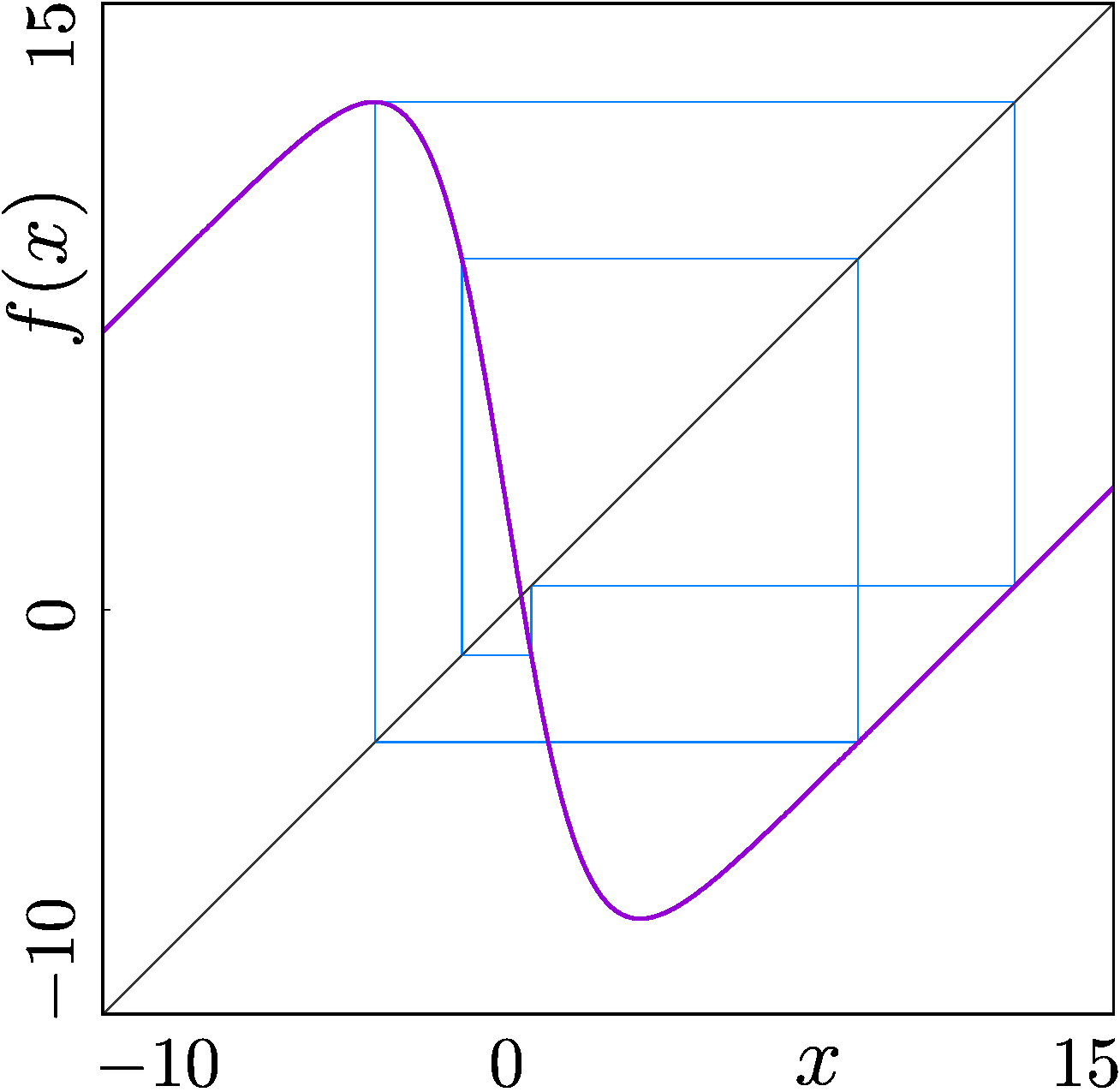}{b}\\[1ex]
  \includesubgraphics{.4\textwidth}{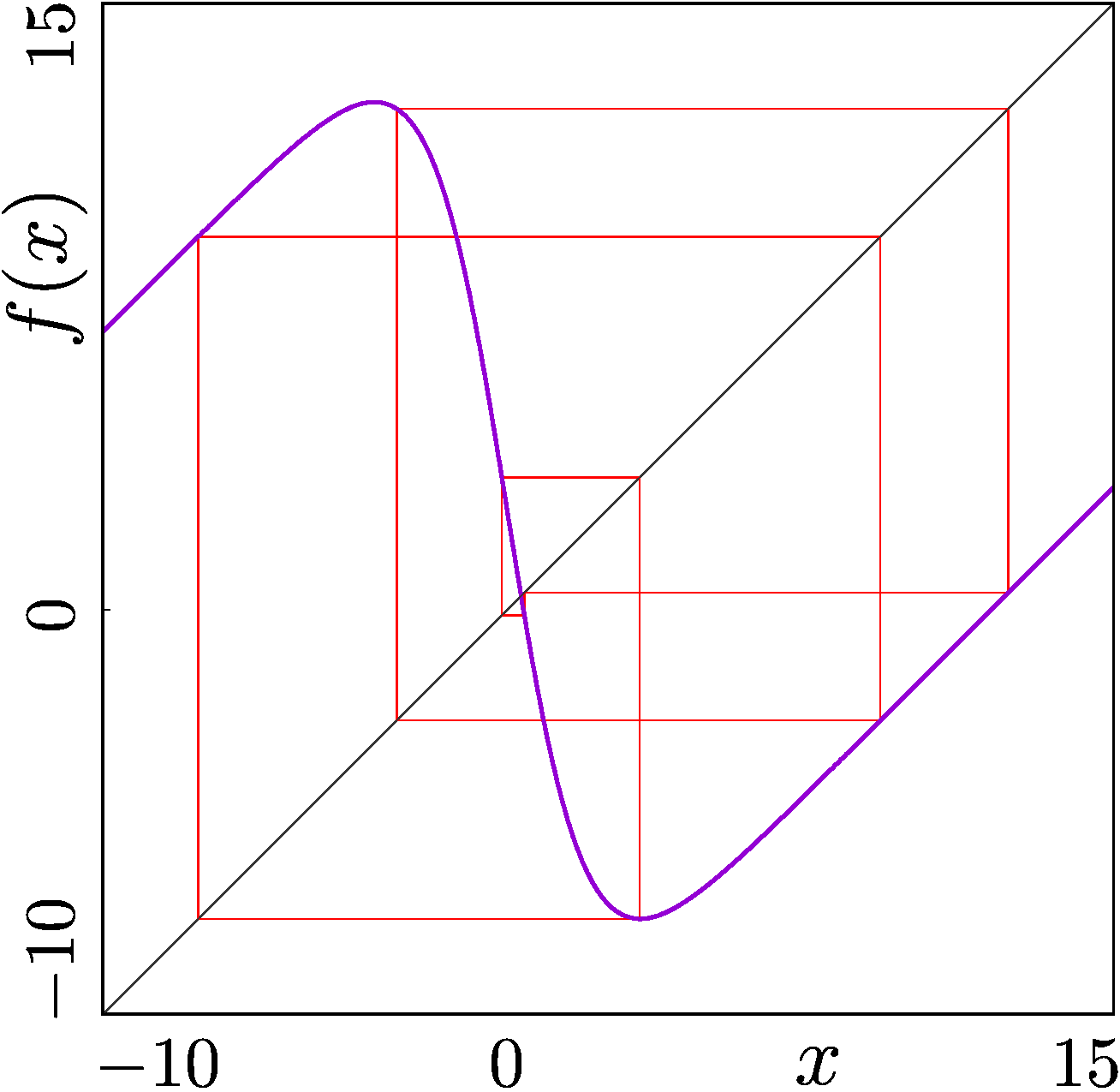}{c}\quad
  \includesubgraphics{.4\textwidth}{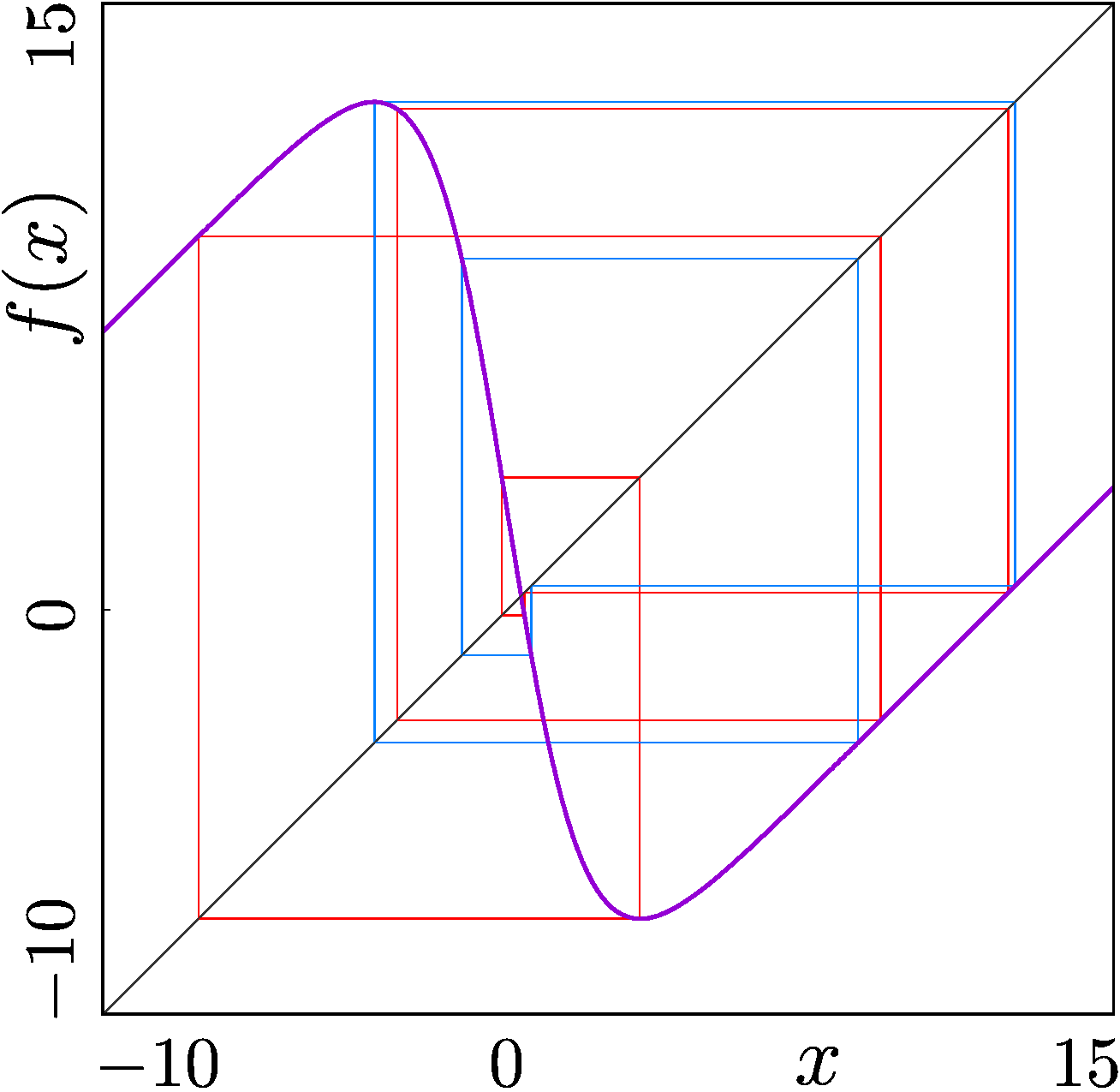}{c}
\end{center}
\caption{\footnotesize Attractors of map~\eqref{glavUrav} existing in the neighborhood of the intersection point of the fold bifurcation curves $\eta^{+1}_5$ and $\eta^{+1}_7$ (point A in Fig.~\ref{fig:BxK:blowup}(c)).
(a) unique broad-band chaotic attractor ($b=-11.9655$, $k=-28.85$);
(b) unique 5-cycle ($b=-11.9642$, $k=-28.8515$);
(c) unique 7-cycle ($b=-11.9672$, $k=-28.853$);
(d) coexisting 5-cycle ($b=-11.9655$, $k=-28.854$);
}
\label{fig:attractors:A}
\end{figure}

As an example, Fig.~\ref{fig:BxK:blowup} presents a few subsequent magnifications of the bifurcation structure shown in Fig.~\ref{bkregpp}. In particular, in Fig.~\ref{fig:BxK:blowup}(c) one can see the bifurcation structure close to the intersection point (marked by A) of the bifurcation curves $\eta^{+1}_5$ and $\eta^{+1}_7$ associated with fold bifurcations of 5- and 7-cycles, respectively.  Examples of attractors of map~\eqref{glavUrav} at parameter values belonging to a neighborhood of point A are shown in Fig.~\ref{fig:attractors:A}.  As one can see, in the quadrant above this point (i.e., before the fold bifurcations occurring at $\eta^{+1}_5$ and $\eta^{+1}_7$), the map has a unique broad-band chaotic attractor (see Fig.~\ref{fig:BxK:blowup}(a)). In the quadrant on the left of this point (i.e., before the fold bifurcation occurring at $\eta^{+1}_5$ but after the fold bifurcation occurring at $\eta^{+1}_7$), the stable 7-cycle is the unique attractor (Fig.~\ref{fig:attractors:A}(b)). Similarly, in the quadrant on the right of point A (i.e., after the fold bifurcation occurring at $\eta^{+1}_5$ but before the fold bifurcation occurring at $\eta^{+1}_7$), the unique attractor of the map is the stable 5-cycle (Fig.~\ref{fig:attractors:A}(c)). As for the quadrant located below the point $A$ (i.e., after both the fold bifurcations occurring at $\eta^{+1}_5$ and $\eta^{+1}_7$), here the stable 7- and 5-cycles coexist (Fig.~\ref{fig:attractors:A}(d)).

\begin{figure}[t!]
\begin{center}
  \includesubgraphics{.4\textwidth}{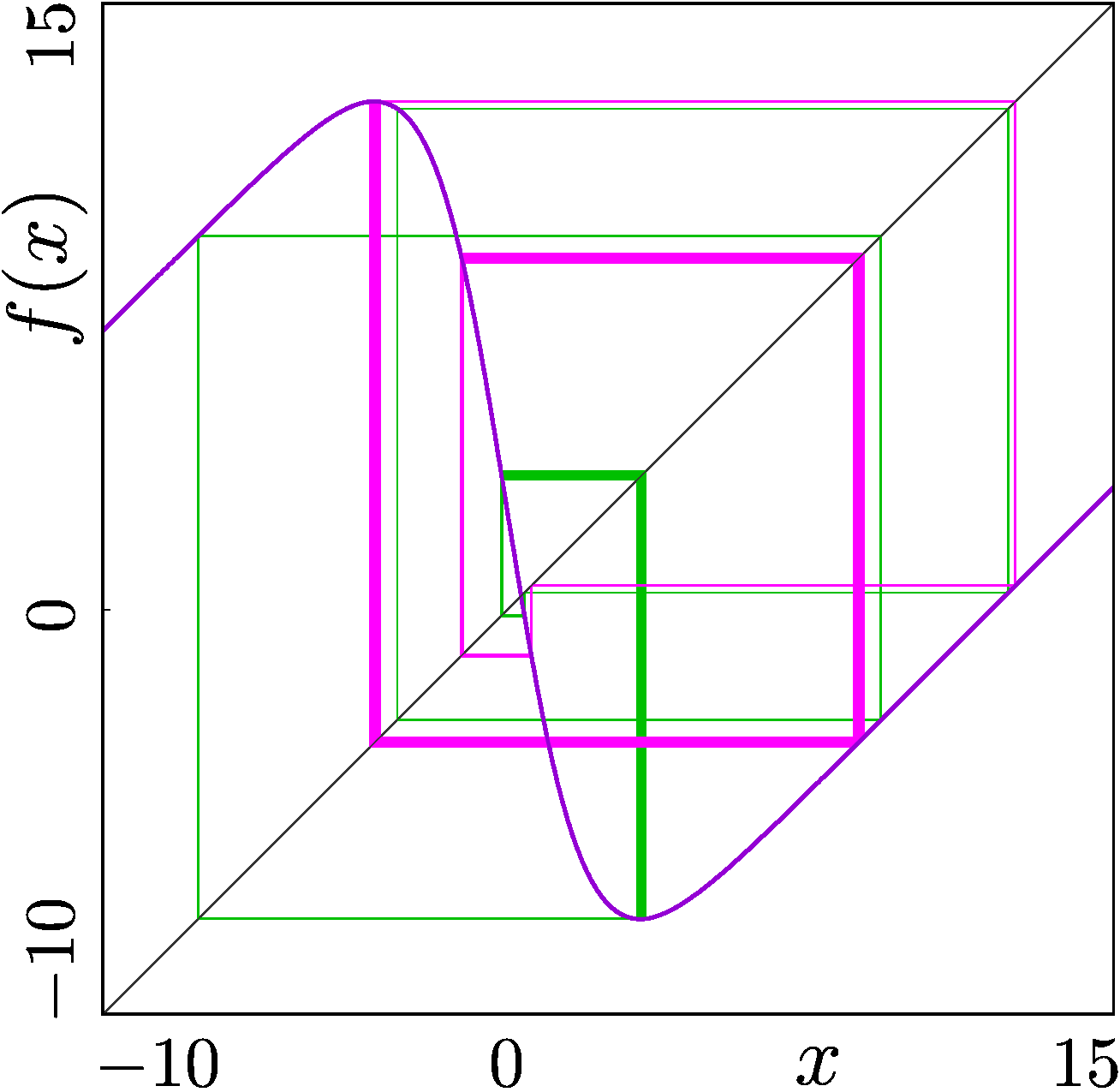}{a}\quad
  \includesubgraphics{.4\textwidth}{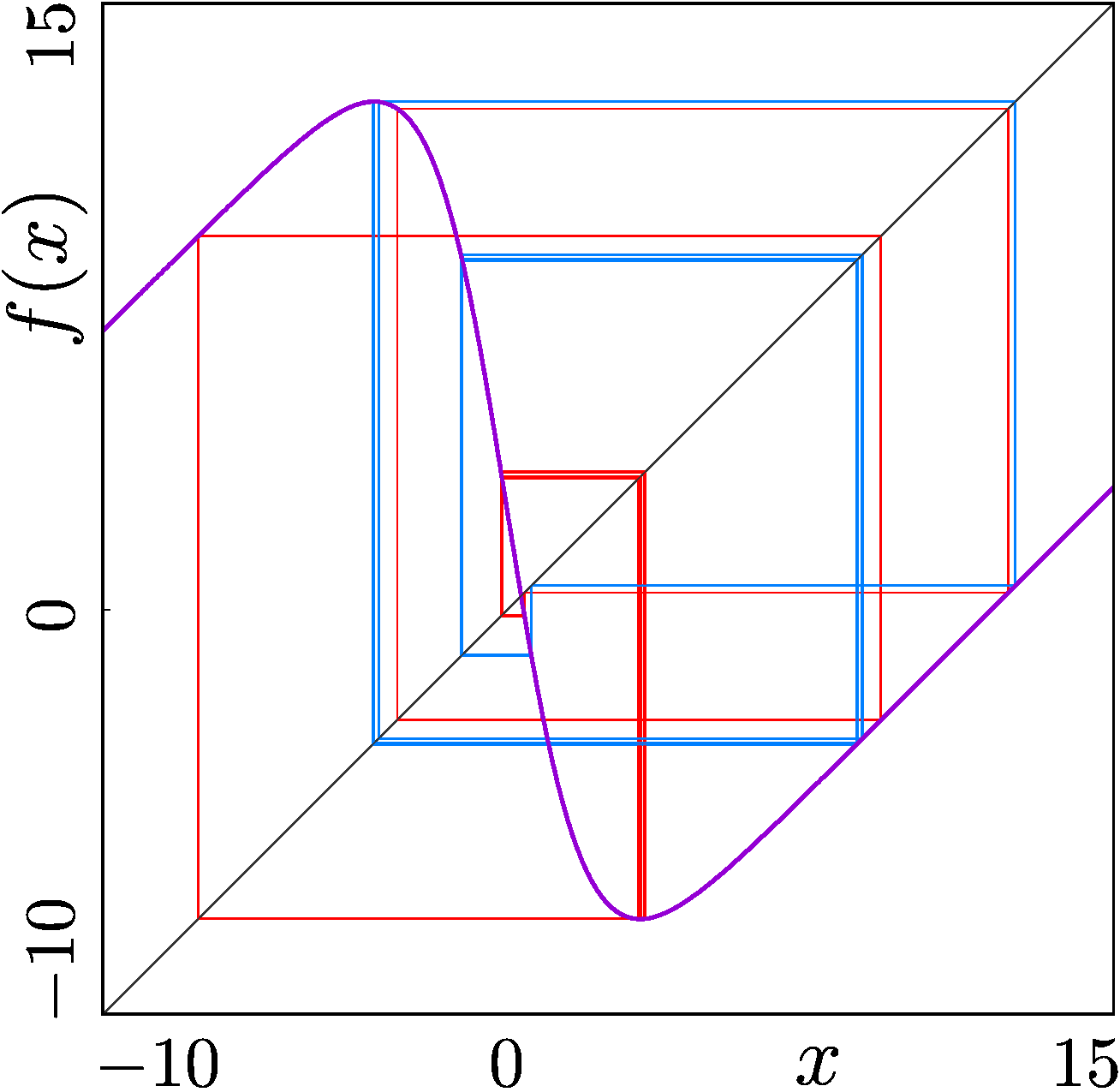}{b}
\end{center}
\caption{\footnotesize
  Coexisting attractors of of map~\eqref{glavUrav}:
  (a) two narrow-band chaotic attractors ($b=-11.9709$, $k=-28.8708$);
  (a) a $20$-cycle and a $56$-cycle ($b=-11.9708$, $k=-28.8695$);
}
\label{fig:attractors:C}
\end{figure}

At the bifurcation curves $\eta^{-1}_5$ and $\eta^{-1}_7$ the stable 5- and 7-cycles appearing at $\eta^{-1}_5$ and $\eta^{-1}_7$, respectively, undergo flip bifurcations.  At the point (marked by B in Fig.~\ref{fig:BxK:blowup}(c)), these curves intersect so that in its neighborhood, one can observe not only the coexistence of stable 5- and 7-cycles (in the quadrant above this point) but also the coexistence of stable 5- and 14- (in the quadrant on the left of this point), 10- and 7- (in the quadrant on the right of this point), and 10- and 14-cycles (in the quadrant below this point).

Under further parameter variation, both period-doubling cascades proceed, followed by complete logistic map scenarios. As a consequence, an arbitrary attractor belonging to one of these scenarios may coexist with an arbitrary attractor belonging to another one. As an example, Fig.~\ref{fig:attractors:C} shows a pair of coexisting narrow-band chaotic attractors and a pair of coexisting cycles of periods 20 and 56.

\section{\label{sec:conclusions}Conclusion}
We considered a model of a two-predators-one-prey system. Previously it was demonstrated that, under certain assumptions, the dynamics system could be represented by a one-dimensional bimodal map. In the present work, we explained the bifurcation structure in the 2D parameter space of this map. We identified the region in the parameter space associated with bounded dynamics. Then, we described the domains in the parameter space associated with different attractors. As we have shown, these sets overlap pairwise, leading to bistability.

\section*{Acknowledgements} Viktor Avrutin was supported by DFG, AV 111/2-2. Sergey Kryzhevich was supported by Gda\'{n}sk University of Technology by the DEC 14/2021/IDUB/I.1 grant under the Nobelium - ‘Excellence Initiative - Research University’ program. Authors dedicate the paper to the memory of Gennadiy Alexeevich Leonov.


\end{document}